\documentclass[11pt]{article}
\usepackage{amsmath,amsfonts, amssymb,amsthm}
\def\R{{\mathbb{R}}}

\def\ep{\varepsilon}
\def\phi{\varphi}

\def\be{\begin{equation}}
\def\en{\end{equation}}
\def\bee{\begin{eqnarray*}}
\def\ene{\end{eqnarray*}}
\def\ben{\begin{eqnarray*}}
\def\een{\end{eqnarray*}}

\title{\protect \large \bf
HYPERBOLIC MEASURES \\ ON INFINITE DIMENSIONAL SPACES
\thanks{Key words: Hyperbolic (convex) measures, dimension, localization, 
dilation of sets}
}

\author {Sergey G. Bobkov  \ 
and James Melbourne}

\begin{document}

\maketitle

\begin{abstract}
\hskip-5.5mm
Localization and dilation procedures are discussed for infinite dimensional 
$\alpha$-concave measures on abstract locally convex spaces
(following Borell's hierarchy of hyperbolic measures).
\end{abstract}


\section{Introduction}
\setcounter{equation}{0}

\noindent
The purpose of this note is to review some results about the localization techniques 
and hyperbolic measures on $\R^n$ and to discuss possible extensions to the setting of 
abstract (infinite dimensional) locally convex spaces. As a starting point, let 
us recall the so-called ``Localization Lemma" which is due to Lov\'asz and Simonovits.

\vskip5mm
{\bf Theorem 1.1} ([L-S]). 
{\it Let $u,v:\R^n \rightarrow \R$ be lower semi-continuous, integrable functions 
such that
\be
\int_{\R^n} u(x)\,dx > 0, \qquad \int_{\R^n} v(x)\,dx > 0.
\en
Then, for some points $a,b \in \R^n$ and a positive affine function $l$ on 
$(0,1)$,
\be
\int_0^1 u((1-t)a + tb)\,l(t)^{n-1}\,dt > 0, \qquad 
\int_0^1 v((1-t)a + tb)\,l(t)^{n-1}\,dt > 0.
\en
}

There are some other variants of this theorem, for example, when the first integrals
involving the function $u$ are vanishing, cf. [K-L-S]. The approach of
Lov\'asz and Simonovits was based on the concept of a needle coming as
result of a localization (or bisection) procedure. Later Fradelizi and Gu\'edon [F-G1]
proposed an alternative geometric argument with involvement of Krein-Milman's theorem, 
cf. also [F-G2].

Theorem 1.1 is a powerful tool towards certain integral relations in $\R^n$; 
it allows reduction to related inequalities in dimension $n=1$. 
It is therefore not surprising that this theorem has 
found numerous applications in different problems of multidimensional Analysis 
and Geometry, such as isoperimetric problems over convex bodies, 
log-concave and more general hyperbolic measures, as well as Khinchine and 
dilation-type inequalities (cf. [K-L-S], [G], [B1,2,3,6], [N-S-V], [B-N], 
[F], [B-M]).
In many such applications, one considers integrals with respect to measures 
that are different than the Lebesgue measure on $\R^n$, and
therefore a more flexible version of Theorem 1.1 involving other measures 
would be desirable. In addition, having in mind dimension free
phenomena and applications to random processes with hyperbolic distributions, 
it is useful to avoid reference to the dimension and 
to obtain similar statements about spaces and measures of an infinite dimension.

A positive Radon measure $\mu$ on a locally convex space $E$ is called 
$\alpha$-concave ($-\infty \leq \alpha \leq \infty$) if, for all 
non-empty Borel sets $A$ and $B$ in $E$ and for all $0<t<1$,
\be
\mu_*\big((1-t)A + tB\big) \geq \Big[(1-t)\mu(A)^\alpha + t\mu(B)^\alpha\Big]^{1/\alpha}.
\en
Here, $(1-t)A + tB = \{(1-t)x + ty: x \in A, \ y \in B\}$ stands for the Minkowski
weighted sum, and $\mu_*$ is the inner measure (for a possible case 
when $(1-t)A + tB$ is not Borel measurable). By the Radon property, (1.1) may 
equivalently be stated for all non-empty compact subsets of $E$,
and then the inner measure is not needed.

Any measure supported on a one-point set is $\infty$-concave.
In all other cases, necessarily $\alpha \leq 1$.
For example, the Lebesgue measure on $\R$ is $1$-concave.
More generally, the Lebesgue measure on $\R^n$ is $\frac{1}{n}$-concave,
which is the content of the Brunn-Minkowski theorem.

Inequality (1.3) strengthens with growing $\alpha$. In the limit
case $\alpha = -\infty$, (1.3) becomes corresponding to
\be
\mu_*\big((1-t)A + tB\big) \geq \min\{\mu(A),\mu(B)\},
\en
which describes the largest class. Such measures $\mu$ are called convex or hyperbolic.
One important case is also $\alpha = 0$, for which (1.3) is understood as
$$
\mu_*\big((1-t)A + tB\big) \geq \mu(A)^{1-t} \mu(B)^t.
$$
Then the measure $\mu$ is called logarithmically concave, or just log-concave.

The class of log-concave measures on $\R^n$ was first considered by 
Pr\'ekopa [Pr] and previously in dimension one by other authors (cf. [I], [D-K-H]). 
The more general classes of $\alpha$-concave measures (in the setting 
of an abstract locally convex space) were introduced by Borell [Bor1]. 
He studied basic properties of $\alpha$-concave measures, including $0-1$ law, 
integrability of norms, convexity properties of measures under convolutions. 
Borell also gave a characterization of the $\alpha$-concavity in terms 
of densities of finite dimensional projections, cf. [Bor1-2], and also [B-L].

Hyperbolic measures are known to satisfy many other important properties that are
usually expressed in terms of various relations such as, for example, Khinchin-type 
inequalities for polynomials of a bounded degree. What is remarkable, 
most of them involve only the convexity parameter $\alpha$ and do not depend 
on the dimension of the underlying space $E$. It is therefore natural 
to state these relations without any restrictions on $E$ where possible and anyway 
wider than in the popular setting of the Euclidean space $E = \R^n$.
For example, returning to Theorem 1.1, it may be complemented with the following.

\vskip5mm
{\bf Theorem 1.2.} {\it Let $\mu$ be a finite $\alpha$-concave measure on a complete
locally convex space $E$, and let $u,v:E \rightarrow \R$ be lower semi-continuous 
$\mu$-integrable functions such that
\be
\int_{E} u\,d\mu > 0, \qquad \int_{E} v\,d\mu > 0.
\en
Then, for some points $a,b \in E$ and some finite $\alpha$-concave measure $\nu$
supported on the segment $\Delta = [a,b]$,
\be
\int_{\Delta} u\,d\nu > 0, \qquad 
\int_{\Delta} v\,d\nu > 0.
\en
}

Note that lower semicontinuous functions are bounded below on any compact set. 
Hence, their integrals over compactly supported finite measures 
such as (1.2) and (1.6) always exist. As an example, the indicator functions 
of open subsets of $E$ are all lower semicontinuous.

The completeness assumption (meaning that every Cauchy net in $E$ is convergent)
is quite natural. It ensures that the closed convex hull of any compact set in $E$ 
is also compact. In that case any finite Radon measure $\mu$ on $E$ has 
a stronger property
\be
\sup\{\mu(K): K \subset E \ {\rm convex \ compact}\} = \mu(E).
\en
This property is crucial in some applications, but without completeness it is not 
true in general. (Its validity remains unclear e.g. for Radon Gaussian 
measures.)

One can also give a geometric variant of Theorem 1.1 together with a finer 
formulation of Theorem 1.2 in terms of extreme points of the set 
$\mathcal{P}_\alpha(u)$ of all $\alpha$-concave probability measures supported 
on a convex compact set $K \subset E$ and such that $\int u\, d\mu \geq 0$ 
(for a continuous function $u$ on $K$). As we already mentioned, 
this interesting approach to localization was developed by 
Fradelizi and Gu\'edon [F-G1]. It was shown there that in case $E = \R^n$ and 
$\alpha \leq \frac{1}{2}$, any extreme point in $\mathcal{P}_\alpha(u)$ is 
either a mass point or it is supported on an interval $\Delta \subset K$ 
with density $l^{(1-\alpha)/\alpha}$ 
(where $l$ is a non-negative affine function on $\Delta$).
As will be explaned in Section 3, this property extends to general 
locally convex spaces, and then it easily implies Theorem 1.2. 

One interesting application of Theorem 1.2 may be stated in terms of the following
operation proposed in [N-S-V]. Given a Borel subset $A$ in a closed convex set 
$F \subset E$ and a number $\delta \in [0,1]$, define
$$
A_\delta = \Big\{x \in A: m_\Delta(A) \geq 1 - \delta \ 
\ {\rm for \ any \ interval} \ \Delta \subset F \ {\rm such \ that} 
\ x \in \Delta\Big\},
$$
where $m_\Delta$ denotes the normalized one-dimensional 
Lebesgue measure on $\Delta$. 

For example, if $F=E$ and $A$ is the complement to a centrally symmetric, open, 
convex set $B \subset E$, then $A_\delta = E \setminus (\frac{2}{\delta} - 1) B$ 
represents the complement to the corresponding dilation of $B$.

\vskip5mm
{\bf Theorem 1.3.} {\it Let $\mu$ be an $\alpha$-concave probability measure 
on a complete locally convex space $E$ supported on a closed convex set $F$ 
$(-\infty < \alpha \leq 1)$. For any Borel set $A$ in $F$ and for all 
$\delta \in [0,1]$ such that $\mu^*(A_\delta)>0$,
\be
\mu(A) \, \geq \, 
\big[\,\delta \mu^*(A_\delta)^\alpha + (1-\delta)\big]^{1/\alpha}.
\en
}

Here $\mu^*$ denotes the outer measure (which is not needed, when $E$ is a Fr\'echet 
space). This relation resembles very much the definition (1.3).
 
In the important particular case $\alpha=0$ (i.e., for log-concave measures), 
(1.8) becomes
$$
\mu(A)\, \geq\, \mu^*(A_\delta)^\delta.
$$
It was discovered by Nazarov, Sodin and Vol'berg [N-S-V]. The extension of this
result to the class of $\alpha$-concave measures in the form (1.8) 
is settled in [B-N] and [F], still for finite dimensional spaces. 
All proofs are essentially based on Theorem 1.1 or its modifications to reduce 
(1.8) to dimension one (although the one dimensional case appears to be 
rather delicate). Here we make another step removing the dimensionality of 
the space assumption, cf. Section 6.

The organization of this note is as follows. In Section 2 we recall basic general
facts about $\alpha$-concave measures, including Borell's characterization
of the $\alpha$-concavity in terms of densities, and describe several examples.
Sections 3-4 are devoted to the extension of Fradelizi-Gu\'edon's theorem and
Lov\'asz-Simonovits' bisection argument. In particular, the existence
of needles which we understand in a somewhat weaker sense is proved for probability 
measures on Fr\'echet spaces that satisfy the zero-one law. This can be used
as an approach towards Theorems 1.1-1.2, but potenitally may have a wider
range of applications. Finally, Sections 5-6 are devoted to Theorem 1.3, which is
then illustrated in the problem of large and small deviations (Section 7).

We do not try to describe in detail results and techniques in dimension one, 
but mainly focus on their extensions to the setting of infinite dimensional spaces.


\vskip2mm
\section{Support, dimension and characterizations}
\setcounter{equation}{0}

\noindent
The support $H_\mu = {\rm supp}(\mu)$ of any Radon measure $\mu$ on $E$ is 
defined as the smallest closed subset of $E$ of full measure, so that
$\mu(E \setminus H_\mu) = 0$. If $\mu$ is hyperbolic, then the set $H_\mu$ 
is necessarily convex, as follows from (1.4). This set has some dimension 
$$
k = {\rm dim}(\mu) = {\rm dim}(H_\mu), 
$$
finite or not, which is called the dimension of the hyperbolic measure $\mu$.
If it is finite, absolute continuity of $\mu$ will always be understood
with respect to the $k$-dimensional Lebesgue measure on $H_\mu$.

First, let us recall an important general property of hyperbolic measures
proven by Borell.

\vskip5mm
{\bf Theorem 2.1} ([Bor1]). {\it If $\mu$ is a hyperbolic probability measure on
a locally convex space $E$, then for any additive subgroup $H$ of $E$, ether
$\mu_*(H) = 0$ or $\mu_*(H) = 1$.
}

\vskip5mm
In particular, any $\mu$-measurable affine subspace of $E$ has measure either
zero or one.

In [Bor1-2], Borell also gave a full description of $\alpha$-concave measures.
Similarly to (1.3), a non-negative function $f$ defined on a convex subset $H$ of 
$E$ is called $\beta$-concave, if it satisfies
\be
f\big((1-t)x + ty\big) \geq \Big[(1-t)f(x)^\beta + tf(y)^\beta\Big]^{1/\beta}
\en
for all $t \in (0,1)$ and all points $x,y \in H$ such that $f(x)>0$ and $f(y)>0$.
The right-hand side is understood in the usual limit sense for the values
$\beta = -\infty$, $\beta = 0$ and $\beta = \infty$.

\vskip5mm
{\bf Theorem 2.2} ([Bor1]).
{\it If $\mu$ is a finite $\alpha$-concave measure on $\R^n$ of dimension 
$k = {\rm dim}(\mu)$, then $\alpha \leq \frac{1}{k}$. 
Moreover, $\mu$ is absolutely continuous with respect to Lebesgue measure on $H_\mu$ 
and has density $f$ which is positive, finite, and 
$\beta$-concave on the relative interior of $H_\mu$, where
$$
\beta = \frac{\alpha}{1 - \alpha k}.
$$
Conversely, if a measure $\mu$ on $\R^n$ is supported on a convex set $H$
of dimension $k$ and has there a positive, $\beta$-concave density $f$ with
$\beta \geq -\frac{1}{k}$, then $\mu$ is $\alpha$-concave.
}

\vskip5mm
Note that $\beta$ is continuously increasing in the range 
$[-\frac{1}{k},\infty]$, when $\alpha$ is varying in $[-\infty,\frac{1}{k}]$. 

In the extremal case $\alpha = \frac{1}{k}$, the density $f(x) = \frac{d\mu(x)}{dx}$ 
is $\infty$-concave and is therefore constant: Up to a factor, 
$\mu$ must be the $k$-dimensional Lebesgue measure on $H_\mu$.

More generally, if $\alpha \leq \frac{1}{k}$, $\alpha \neq 0$, the density has 
the form
$$
f(x) = V(x)^{\frac{1}{\alpha} - k}
$$
for some function $V:\Omega \rightarrow (0,\infty)$ on the relative interior
$\Omega$ of $H_\mu$, which is concave in case $\alpha > 0$, and is convex 
in case $\alpha < 0$. In particular, the formula 
$$
f(x) = V(x)^{-k}
$$
describes all $k$-dimensional hyperbolic measures $(\alpha = -\infty)$. 
If $\alpha = 0$, then necessarily
$
f(x) = e^{-V(x)}
$
for some convex function $V:\Omega \rightarrow \R$.

As for general locally convex spaces, another theorem due to Borell
reduces the question to Theorem 2.2.

\vskip5mm
{\bf Theorem 2.3} ([Bor1]).
{\it A Radon probability measure $\mu$ on the locally convex space $E$ 
is $\alpha$-concave, if and only if the image of $\mu$ under any linear continuous 
map $T:E \rightarrow \R^n$ is an $\alpha$-concave measure on $\R^n$.
}

\vskip5mm
For special spaces in this characterization one may consider linear continuous maps 
$T$ from a sufficiently rich family. For example, when $E = C[0,1]$ is the
Banach space of all continuous functions on $[0,1]$ with the maximum-norm, 
the measure $\mu$ is $\alpha$-concave, if and only if the image of $\mu$ under 
any map of the form
$$
Tx = (x(t_1),\dots,x(t_n)), \qquad x \in C[0,1], \ \ t_1,\dots,t_n \in [0,1],
$$
is an $\alpha$-concave measure on $\R^n$. Simialrly, when $E = \R^\infty$
is the space of all sequences of real numbers (with the product topology),
it is sufficient to consider the standard projections
\be
T_nx = (x_1,\dots,x_n), \qquad x = (x_1,\dots,x_n,\dots) \in \R^\infty.
\en

\vskip2mm
The next general observation emphasizes that infinite dimensional
$\alpha$-concave measures may not have a positive parameter of convexity.
Apparently, it was not stated explicitly in the literature, so we include the proof.
As usual, $E'$ denotes the dual spaces of all linear continuous functionals on $E$.

\vskip5mm
{\bf Theorem 2.4.}
{\it For $\alpha > 0$, any $\alpha$-concave finite measure $\mu$ on a locally
convex space $E$ has finite dimension and is compactly supported.
}

\vskip5mm
{\bf Proof.} 
First, suppose to the contrary that $\mu$ is infinite dimensional.
We may assume that $H_\mu = {\rm supp}(\mu)$ contains the origin. 
Since $H_\mu$ is not contained in any finite dimensional subspace of $E$, for each 
$n$, one can find linearly independent 
vectors $v_1,\dots,v_n \in H_\mu$. Each point $x \in E$ has a representation 
$x = c_1(x) v_1 + \dots + c_n(x) v_n + y$ with some $c_i \in E'$, where 
$y = y(x)$ is linearly independent of all $v_i$ (cf. [R], Lemma 4.21).
Consider the linear map $T(x) = (c_1(x),\dots,c_n(x))$, which is continuously 
acting from $E$ to $\R^n$. Then the image $\nu = \mu T^{-1}$ of $\mu$ 
is a finite $\alpha$-concave measure on $\R^n$.

Let us see that $\nu$ is full dimensional. Otherwise, $\nu$ is supported on 
some hyperplane in $\R^n$ described by the equation
$a_1 y_1 + \dots + a_n y_n = a_0$, where the coefficients $a_i \in \R$
are not all zero. Moreover, since $0 \in H_\mu$, any neighborhood
of $0$ has a positive $\mu$-measure, so 
$$
\mu\{x \in E: |T(x)| < \ep\} > 0,
$$
for any $\ep>0$. Hence, necessarily $a_0 = 0$. This implies that $\mu$ is 
supported on the closed linear subspace $H$ of $E$ described by the equation
$a_1 c_1(x) + \dots + a_n c_n(x) = 0$. Here, at least one of the coefficient, 
say $a_i$, is non-zero. Since $c_i(v_i) = 1 \neq 0$, we obtain that
$v_i \notin H$. But this would mean that $H_\mu \cap H$ is a proper closed subset 
of the support of $\mu$, while $H_\mu$ has a full $\mu$-measure, a contradiction.

Hence, ${\rm dim}(\nu) = n$. By Theorem 2.2, this gives 
$\alpha \leq \frac{1}{n}$, and since $n$ was arbitrary, we conclude that 
$\alpha \leq 0$ which contradicts to the hypothesis $\alpha>0$.

Thus, $\mu$ must be supported on a finite dimensional affine subspace 
$H \subset E$. To prove compactness of the support, we may assume that 
$H = E = \R^n$ and ${\rm dim}(\mu) = n$. Then, $\mu$ is supported on an open 
convex set $\Omega \subset \R^n$, where it has density of the form
$$
f(x) = V(x)^\gamma, \qquad \gamma = \frac{1}{\alpha} - n \quad
\big(\,0 < \alpha \leq \frac{1}{n}\,\big),
$$
for some concave function $V:\Omega \rightarrow (0,\infty)$. 
The case $\gamma = 0$ is possible, but then $f(x) = c$ for some constant $c>0$,
which implies $\mu(\R^n) = c\, |\Omega|$. Since $\mu$ is finite, $\Omega$ has 
to be bounded, and so $H_\mu = {\rm clos}(\Omega)$ is compact. 

Now, let $\gamma>0$. Suppose that $\Omega$ is unbounded (to justify
several notations below).
It is known (cf. e.g. [B5]) that $f(x) \rightarrow 0$ uniformly as 
$|x| \rightarrow \infty$ ($x \in \Omega$). In particular, $f$ is bounded, 
that is, $A = \sup_{x \in \Omega} f(x)$ is finite. Here, we may assume 
that the sup is asymptotically attained at $x_0 = 0$ for some sequence 
$x_l \rightarrow 0$, 
$x_l \in \Omega$. Choose $r>0$ so that $f(x) < \frac{1}{2}\, A$ or 
$V(x) < (\frac{1}{2}\, A)^{1/\gamma}$ whenever $|x| \geq r$ and $x \in \Omega$. 
For such $x$, the sequence 
$$
\lambda_l(x) = \sup\{\lambda > 1: x_l + \lambda (x-x_l) \in \Omega\}
$$
has a limit $\lambda(x) = \sup\{\lambda: \lambda_l \in \Omega\} > 1$.
Consider the convex functions
$$
\psi_l(\lambda) =  V(x_l) - V(x_l + \lambda (x-x_l)), \qquad 
0 \leq \lambda < \lambda_l(x).
$$ 
We have $\psi_l(0)=0$ and 
$\psi_l(1) = V(x_l) - V(x) > C = A^{1/\gamma} (1 - 2^{-1/\gamma})$, for all $k$ 
large enough. Hence, $\psi_l(\lambda) \geq C \lambda$, for all admissible 
$\lambda \geq 1$, and letting $l \rightarrow \infty$, we obtain
$$
C \lambda \leq A^{1/\gamma} - V(\lambda x), \qquad 
1 \leq \lambda < \lambda(x) \ \ (|x| \geq r, \ x \in \Omega).
$$
But $V$ is non-negative, so necessarily 
$\lambda(x) \leq \frac{1}{1 - 2^{-1/\gamma}}$. This proves boundedness of $\Omega$.
\qed

\vskip5mm
{\bf Examples 2.5.}

1. The normalized Lebesgue measure on every convex body $K \subset \R^n$ is 
$\frac{1}{n}$-concave.

2. Any Gaussian measure on a locally convex space $E$ is log-concave. 
In particular, the Wiener measure on $C[0,1]$ is such.

3. The standard Cauchy measure $\mu_1$ on $\R$ with density 
$f(x) = \frac{1}{\pi (1 + x^2)}$ is $\alpha$-concave with $\alpha = -1$ 
(which is optimal). More generally, the $n$-dimensional Cauchy measure $\mu_n$ on 
$\R^n$ with density 
$$
f_n(x) = \frac{c_n}{(1 + |x|^2)^{(n+1)/2}}
$$
is $(-1)$-concave ($c_n$ is a normalizing constant so that $\mu_n$ is probability).

4. Although the above density $f_n$ essentially depends on the dimension,
the measure $\mu_n$ has a dimension-free essense. All marginals of $\mu_n$ coincide
with $\mu_1$ and moreover, there is a unique Borel probability measure $\mu$
on $\R^\infty$ (an infinite dimensional Cauchy measure) which is pushed forward to 
$\mu_n$ by the standard projection $T_n$ from (2.2). This measure can also be
introduced as the distribution of the random sequence
$$
X = \Big(\frac{Z_1}{\zeta},\frac{Z_2}{\zeta},\dots\Big),
$$
where the random variables $\zeta,Z_1,Z_2,\dots$ are independent and all have a standard
normal distribution. Thus, $\mu$ is $(-1)$-concave on $\R^\infty$.

5. This example is mentioned in [Bor1]. Given $d>0$ (real), let $\chi_d$ be a positive 
random variable such that $\chi_d^2$ has the $\chi^2$-distribution with
$d$ degrees of freedom, i.e., with density 
$$
f_d(r) = \frac{1}{2^{d/2} \Gamma(d/2)}\,r^{d/2 - 1}\,e^{-r/2}, \qquad r > 0.
$$
Let $W$ be the standard Wiener process (independent of $\chi_d$) viewed as a random 
function in $C[0,1]$. Then the random function
$$
X(t) = \frac{\sqrt{d}}{\chi_d}\,W(t), \quad t \in [0,1],
$$
has the distribution $\mu$ which is $\alpha$-concave on $C[0,1]$ with 
$\alpha = -\frac{1}{d}$. It is called the Student measure (and also Cauchy in case $d=1$
similarly to the previous example).


\vskip2mm
\section{Extreme $\alpha$-concave measures}
\setcounter{equation}{0}

\vskip2mm
\noindent
Given a convex compact set $K$ in a locally convex space $E$, denote by
$\mathcal{M}_\alpha(K)$ the collection of all $\alpha$-concave probability 
measures with support contained in $K$. For a continuous function $u$ on $K$,
we consider the subcollection
$$
\mathcal{P}_\alpha(u) = \left\{\mu \in \mathcal{M}_\alpha(K): \int u\, d\mu \geq 0\right\}
$$
together with its closed convex hull $\mathcal{\widetilde P}_\alpha(u)$ 
in the locally convex space $\mathcal{M}(K)$ of all signed Radon measures on $K$ 
endowed with the topology of weak convergence. The latter space is
dual to the space $C(K)$ of all continuous functions on $K$, and
$\mathcal{\widetilde P}_\alpha(u)$ is a convex compact subset of $\mathcal{M}(K)$.

What are extreme points of $\mathcal{\widetilde P}_\alpha(u)$?
Using a general theorem due to D. P. Milman, one can only say that all such 
points lie in $\mathcal{P}_\alpha(u)$ (cf. [B-S-S], p.124, or [Ph] 
for a detail discussion of Krein-Milman's theorem). A full answer to this question 
is given in Fradelizi-Gu\'edon's theorem, which we formulate below in the setting 
of abstract locally convex spaces.

\vskip5mm
{\bf Theorem 3.1.} {\it Given a continuous function $u$ on $K$ and 
$-\infty \leq \alpha \leq 1$, any extreme point $\mu$ in 
$\mathcal{\widetilde P}_\alpha(u)$ has the dimension ${\rm dim}(\mu) \leq 1$. 
Moreover, in case $\alpha \leq \frac{1}{2}$, 

\vskip3mm
$1)$\, $\mu$ is either a mass point at $x \in K$ such that 
$u(x) \geq 0$; or 

$2)$\, $\mu$ is supported on an interval $\Delta = [a,b] \subset K$ with density 
\be
\frac{d\mu(x)}{dm_\Delta(x)} = l(x)^{(1-\alpha)/\alpha}
\en
with respect to the uniform measure $m_\Delta$,
where $l$ is a non-negative affine function on $\Delta$ such that
$\int_a^x u\,d\mu > 0$ and $\int_x^b u\,d\mu > 0$, for all $x \in (a,b)$.
}

\vskip5mm
In particular, any $\alpha$-concave probability measure supported on $K$, belongs 
to the closed convex hull of the family of all one-dimensional $\alpha$-concave 
probability measures supported on $K$ having density of the form (3.1).

We only consider the first assertion of the theorem. The second part is a purely 
one dimensional statement, and we refer to [F-G1].

\vskip5mm
{\bf Proof.} Suppose that a measure $\mu \in \mathcal{P}_\alpha(u)$ has 
the dimension ${\rm dim}(\mu) \geq 2$. For simplicity, let the origin belong 
to the relative interior $G$ of the support $H_\mu$ of $\mu$. Then one may find 
linearly independent vectors $x$ and $y$ such that $\pm x$ and $\pm y$ are all
in $G$. On the linear hull $L(x,y)$ of $x$ and $y$ 
(which is a 2-dimensional linear subspace of $E$), define linear
functionals $\lambda_x$ and $\lambda_y$ by putting
$$
\lambda_x(x)= \lambda_y (y) = 1,
$$
$$
\lambda_x(y)= \lambda_y(x) = 0.
$$
They are continuous, so by the Hahn-Banach theorem, these functionals may be 
extended from $L(x,y)$ to the whole space $E$ keeping linearity and continuity. 

With these extended functionals, we can associate 
$\Lambda_\theta = \theta_1 \lambda_x + \theta_2 \lambda_y$, where
$\theta = (\theta_1, \theta_2) \in \mathbb{S}^1$ (vectors on the unit sphere of $\R^2$).
Note that these functionals are uniformly bounded on $K$, i.e.,
\be
\sup_{\theta} \, \sup_{z \in K}\, |\Lambda_\theta(z)| \, \leq \,
\sup_{z \in K}\, |\lambda_x(z)| + \sup_{z \in K}\, |\lambda_y(z)| \, < \, \infty.
\en
 
Now, following in essense an argument of [F-G1], define the map 
$\Phi: \mathbb{S}^1 \to \mathbb{R}$ by
$$
\Phi(\theta) = \int_{\{\Lambda_\theta \geq 0\}} u\, d\mu.
$$
By the construction, the set $\{\Lambda_\theta = 0\} \cap H_\mu$ represents
a proper closed affine subspace of $H_\mu$. So, $\mu\{\Lambda_\theta = 0\} = 0$
according to Theorem 3.1 (the zero-one law for hyperbolic measures). Hence, 
using (3.2), we may conclude that the map $\Phi$ is continuous. 

In addition, we have the identity $\Phi(\theta)+\Phi(-\theta) = \int u\, d\mu$. 
Hence, the intermediate value theorem implies that there exists $\theta$ 
such that with 
$H_\theta^+ = \{ \Lambda_\theta \geq 0 \}$ and 
$H_\theta^- = \{ \Lambda_\theta \leq 0 \}$, we have
$$
\int_{H_\theta^+} u\,d\mu = \int_{H_\theta^-} u\,d\mu = \frac{1}{2} \int_{E} u\, d\mu.
$$
Necessarily, $t = \mu(H_\theta^-) > 0$ and $\mu(H_\theta^+) > 0$.
Defining $\alpha$-concave probability measures 
$$
\mu_0(A) = \frac{\mu(A \cap H_\theta^+)}{\mu(H_\theta^+)}, \qquad
\mu_1(A) = \frac{\mu(A \cap H_\theta^-)}{\mu(H_\theta^-)},
$$
we arrive at the representation $\mu = (1-t)\mu_0 + t\mu_1$ which means that $\mu$ 
is not extreme. 
\qed

One can now return to Theorem 1.2.

\vskip5mm
{\bf Proof of Theorem 1.2.} Due to the property (1.7), and by the assumption (1.5),
$$
\int_K \min(u,c)\,d\mu > 0, \qquad \int_K \min(v,c)\,d\mu > 0,
$$
for some convex compact set $K \subset E$ and a constant $c>0$. Moreover, since 
the function $\min(u,c)$ is lower semicontinuous and bounded, while $\mu$ is Radon,
$$
\int_K \min(u,c)\,d\mu \, = \, \sup_g \int g\,d\mu,
$$
where the sup is taken over all continuous functions on $K$ such that
$g \leq \min(u,c)$ (cf. e.g. [M], Chapter 2, or [Bog], Chapter 7). A similar 
identity also holds for $\min(v,c)$. This allows us to reduce the statement of 
the theorem to the case where both $u$ and $v$ are continuous on $K$.

In the latter case, let $u_0 = u - \int_K u\,d\mu$.
Consider the functional $T(\mu) = \int_K v\,d\mu$. It is
linear and continuous on $\mathcal{M}(K)$, and therefore being restricted to
$\mathcal{P}_\alpha(u_0)$ it attains maximum at one of the extreme points $\nu$.
Since $\mu \in \mathcal{P}_\alpha(u_0)$, we conclude that
$$
\int_K u_0\,d\nu \geq 0, \quad T(\nu) \geq T(\mu),
$$
so, $\int_K u\,d\nu > 0$ and $\int_K v\,d\nu > 0$ which is (1.6). 
It remains to apply Theorem 3.1.
\qed

\vskip5mm
A similar argument, based also on the second part of Theorem 3.1, yields
Theorem 1.1. Indeed, the $n$-dimensional integrals (1.1) can be restricted 
to a sufficently large closed ball $K \subset \R^n$. The normalized Lebesgue 
measure on $K$ is $\alpha$-concave with $\alpha = \frac{1}{n}$.
Hence, the extreme points in $\mathcal{P}_\alpha(u)$
are at most one dimensional and have densities of the form $l^{n-1}$
(if they are not Dirac measures).


\vskip2mm
\section{Bisection and needles on Fr\'echet spaces}
\setcounter{equation}{0}

\vskip2mm
\noindent
The notion of a needle was proposed by Lov\'asz and Simonovits for the proof 
of Theorem 1.1 (Localization Lemma, cf. also [K-L-S]). Previously, it appeared
implicitly in [P-W] and may be viewed as development of the Hadwiger-Ohmann bisection
approach to the Brunn-Minkowski inequality ([H-O], [B-Z], cf. also [G-M] 
for closely related ideas).

As shown in [L-S], starting from (1.1), one can construct a decreasing sequence 
of compact convex bodies $K_l$ in $\R^n$ that are shrinking to some segment 
$\Delta = [a,b]$ and are such that, for each $l$,
$$
\int_{K_l} u(x)\, dx >0, \qquad \int_{K_l} v(x)\, dx > 0.
$$
Moreover, choosing a further subsequence (if necessary) and applying the 
Brunn-Min\-kow\-s\-ki inequality in $\R^n$, one gets in the limit
\bee
\lim_{l \rightarrow \infty}\, \frac{1}{|K_l|} \int_{K_l} u(x)\, dx
 & = &
\int_{\Delta} \psi^{n-1}(x)\,u(x)\, dx, \\
\lim_{l \rightarrow \infty}\, \frac{1}{|K_l|} \int_{K_l} v(x)\, dx
 & = &
\int_{\Delta} \psi^{n-1}(x)\,v(x)\, dx,
\ene
for some non-negative concave function $\psi$ on $\Delta$. Here $|K_l|$ 
denotes the $n$-dimensional volume, while the integration on the right-hand 
side is with respect to the linear Lebesgue measure on the segment. 
In this way, one may obtain a slightly weaker variant of (1.2) with $\psi$ 
in place of $l$, and with non-strict inequalities. An additional argument 
of a similar flavour was then developed in [L-S] to make $\psi$ affine 
(while the strict inequalities in (1.2) are easily achieved by applying
the conclusion to functions $u - \ep w$ and $v - \ep w$, where $w > 0$ 
is integrable, continuous, and $\ep>0$ is small enough).
The last step shows that for $K_l$ one may take infinitesimal truncated cylinders
with main axis $\Delta$; it is in this sense the limit one dimensional measure 
$l^{n-1}(x)\,dx$ on $\Delta$ may be considered a needle.

The aim of this section is to extend this construction to the setting of
Fr\'echet, i.e., complete metrizable locally convex spaces.
For example, $E$ may be a Banach space, but there also other important spaces
that are not Banach, such as the space $E = \R^\infty$. 
Note that any finite Borel measure on a Fr\'echet space is Radon.

While one cannot speak about the Lebesgue measure when $E$ is infinite dimensional, 
the main hypothesis (1.1) may readily be stated like (1.5) with integration 
with respect to a given (finite) Borel measure $\mu$ on $E$.

The space of all finite Borel measures on $E$ is endowed with the topology of
weak convergence. In particular, $\mu_l \rightarrow \mu$ (weakly), if and only if
$$
\int u\,d\mu_l \rightarrow \int u\,d\mu \qquad ({\rm as} \ \ l \rightarrow \infty)
$$
for any bounded continuous functions $u$ on $E$.
As was noticed in [Bor1], the class of all $\alpha$-concave 
probability measures on $E$ is closed in the weak topology.

\vskip5mm
{\bf Definition 4.1.} Let $\mu$ be a finite Borel measure on $E$. 
A Borel probability measure $\nu$ will be called a needle of $\mu$, if it is 
supported on a segment $[a,b] \subset E$ and can be obtained as the weak limit
of probability measures
$$
\mu_l(A) = \frac{1}{\mu(K_l)}\,\mu(A \cap K_l), \qquad (A \ {\rm is \ Borel}),
$$
where $K_l$ is some decreasing sequence of convex compact sets in $E$ of positive
$\mu$-measure such that $\cap_l K_l = [a,b]$.

\vskip5mm
Here, all $\mu_l$ represent normalized restrictions of $\mu$ to $K_l$.
In particular, all needles of a given $\alpha$-concave measure are $\alpha$-concave,
as well. We do not require that $K_l$ be asymptotically close to 
infinitesimal truncated cylinders.

\vskip5mm
{\bf Definition 4.2.}
One says that a Borel probability measure $\mu$ on $E$ satisfies the zero-one law,
if any $\mu$-measurable affine subspace of $E$ has $\mu$-measure either 0 or 1.

\vskip5mm
For example, this important property holds true for all (Radon) Gaussian measures.
More generally, it is satisfied by any hyperbolic probability measure, as follows
from Borell's Theorem 2.1. 

With these definitions, Theorem 1.2 admits the following refinement.

\vskip5mm
{\bf Theorem 4.3.} {\it Suppose that a Borel probability measure $\mu$ on 
a Fr\'echet space~$E$ satisfies the zero-one law. Let $u,v:E \rightarrow \R$ 
be lower semi-continuous $\mu$-integrable functions such that
$$
\int u\,d\mu > 0, \qquad \int v\,d\mu > 0.
$$
Then, these inequalities also hold for some needle $\nu$ of $\mu$.
Moreover, if $\mu$ is supported on a closed convex set $F$, then $\nu$ may 
be chosen to be supported on $F$, as well.
}

\vskip5mm
First assume that $E$ is a separable Banach space with norm $\|\cdot\|$, 
and let $E'$ denote the dual space (of all linear continous functionals on $E$)
with norm $\|\cdot\|_*$. Suppose that any proper closed affine subspace of $E$ 
has $\mu$-measure zero.
In this case, for the proof of Theorem 4.3 we use the construction similar
to the one from the proof of Theorem 3.1.

Given 3 affinely independent points $x,y,z$ in $E$, 
define linear functionals $\lambda_x$ and $\lambda_y$ on the linear hull 
$L_z(x,y)$ of $x-z$ and $y-z$ (which is a 2-dimensional linear subspace of $E$), 
by putting
\be
\lambda_x(x-z)= \lambda_y (y-z) = 1,
\en
\be
\lambda_x(y-z)= \lambda_y(x-z) = 0.
\en
By the Hahn-Banach theorem, these functionals may be extended by linearity to the 
whole space $E$ without increasing their norms. This will always be assumed below.

\vskip5mm
{\bf Lemma 4.4.} {\it Let $\{(x_n,y_n,z_n)\}_{n \geq 1}$ be
affinely independent points in the Banach space $E$ such that $x_n \rightarrow x$, 
$y_n \rightarrow y$, $z_n \rightarrow z$, where $x,y,z$ are also 
affinely independent. Then the corresponding linear functionals $\lambda_{x_n}$ and
$\lambda_{y_n}$ have uniformly bounded norms, i.e.,
$$
\sup_{n \geq 1} \|\lambda_{x_n}\|_* < \infty, \qquad
\sup_{n \geq 1} \|\lambda_{y_n}\|_* < \infty.
$$
}

\vskip5mm
{\bf Proof.} Define the lines 
\bee
L_z(x) & = & \{z + r(x-z): r \in \mathbb{R}\}, \\
L_z(y) & = & \{z + r(y-z): r \in \mathbb{R}\}.
\ene
Then, for $w \in L_z(x,y)$, $\|w\| \leq 1$,
$$
|\lambda_x(w)| \leq {\rm dist}^{-1}(x,L_z(y)), \quad
|\lambda_y(w)| \leq {\rm dist}^{-1}(y,L_z(x)),
$$
where we use the notation ${\rm dist}(w,A) = \inf\{\|w-a\|: a \in A\}$ (the shortest 
distance from a point to the set). The extended linear functionals should thus 
satisfy the above inequalities on the whole space $E$ for all $\|w\| \leq 1$, i.e.,
\be
\|\lambda_x\|_* \leq {\rm dist}^{-1}(x,L_z(y)), \quad
\|\lambda_y\|_* \leq {\rm dist}^{-1}(y,L_z(x)).
\en

Next, by shifting, one may assume that $z=0$, in which case $x$ and $y$ are 
linearly independent and in particular $\|x\|>0$ and $\|y\| > 0$.
Using (4.3), it is enough to show that 
$$
{\rm dist}(x_n,L_{z_n}(y_n))\geq c, \quad {\rm for \ all} \ n \geq n_0, 
$$ 
with some $n_0$ and $c>0$. Indeed, take an arbitrary point $w = z_n + r(y_n - z_n)$ in 
$L_{z_n}(y_n)$, $r \in \R$. By the triangle inequality,
$$
\|x_n - w\| \geq |r| \, \|y_n - z_n\| - \|x_n - z_n\| \geq 2\|x_n - z_n\|,
$$
where the last inequality holds whenever 
$|r| \geq 3\,\frac{\|x_n - z_n\|}{\|y_n - z_n\|}$. Hence, by the convergence assumption,
$$
\|x_n - w\| \geq \|x\|, \qquad {\rm for} \ \ 
|r| \geq r_0 = 4\,\frac{\|x\|}{\|y\|}, \ n \geq n_0.
$$
In case $|r| \leq r_0$, again by the triangle inequality,
\bee
\|x_n - w\| 
 & \geq &
\|x - (z + ry)\| - \|x_n - x\| - |r| \, \|y_n - y\| - r \|z_n\| \\
 & \geq &
{\rm dist}(x,L_{z}(y)) - \|x_n - x\| - r_0 \, \|y_n - y\| - r_0 \|z_n\|.
\ene
Here, the right-hand side is also separated from zero for sufficiently large $n$.

By a similar argument, ${\rm dist}(y_n,L_{z_n}(x_n))\geq c$, for all $n \geq n_0$.
\qed

\vskip5mm
{\bf Proof of Theorem 4.3.} We begin with a series of reductions assuming without 
loss of generality that $\mu(E) = 1$. 

Any Fr\'echet space with Radon probability measure $\mu$ has a subspace $E_0$ 
such that $\mu(E_0)=1$, and in addition there exists a norm $\| \cdot \|$ on 
$E_0$ with respect to which $E_0$ is a separable reflexive Banach space 
whose closed balls are compact in $E$ (see [Bog], Theorem 7.12.4).  

In particular, all Borel subsets of $E_0$ are Borel in $E$.
By the zero-one law (turning to a smaller subspace if necessary), we may assume 
that any proper affine subspace of $E_0$ which is closed for the topology of $E_0$
has measure zero. That is, for all $l \in E_0'$,
\be
\mu\{l = c\} = 0, \quad c \in \R.
\en 

Second, it suffices to assume that the support of $\mu$ is compact, metrizable and 
convex. Indeed, by Ulam's theorem, there is an increasing sequence of compact 
sets $K_n \subset E_0$ such that $\mu(\cup_n K_n) = 1$.
The closed convex hull of any compact set in $E_0$ is compact
(which is true in any Banach and more generally complete locally convex spaces, cf.
e.g. [K-A]). Therefore, all $K_n$ may additionally be assumed to be convex.
These sets will also be compact in $E$, so that the associated weak topologies
in the spaces of Borel probability measures on $K_n$ coincide, as well.
By the dominated convergence theorem,
$$
\lim_{n \to \infty} \int_{K_n} u\, d\mu = \int_{E} u\, d\mu, \qquad 
\lim_{n \to \infty} \int_{K_n} v\, d\mu = \int_{E} v\, d\mu,
$$
so that $\int_{K_n} u\, d\mu > 0$ and $\int_{K_n} v\, d\mu > 0$ for large $n$. 
Hence, an application of the theorem to $\mu$ restricted and normalized to $K_n$ 
would provide the desired one dimensional measure $\nu$, a needle of $\mu_n$ and 
therefore of $\mu$ itself.

Thus, from now on, we may assume that $E$ is a separable Banach space, and
$\mu$ is a Borel probability measure on $E$ which is supported on a convex compact
set $K \subset E$ and is such that (4.4) holds true for all $l \in E'$.

We need only to prove the existence of $\nu$ such that 
$\int u\, d\nu \geq 0$ and $\int v\, d\nu \geq 0$.
Since in this case we may apply the superficially weaker result to $u - \ep$ and 
$v - \ep$ for an $\ep >0$ chosen small enough to preserve the hypothesis.

In addition, it suffices to prove the result when $u$ and $v$ are both continuous.  
To see this, take $u_n$ and $v_n$ to be sequences continous functions increasing 
to lower semicontinuous $u$ and $v$ respectively. By the monotone convergence, 
$\lim_{n \to \infty} \int u_n\, d\mu = \int u\, d\mu > 0$ and 
$\lim_{n \to \infty} \int v_n\, d\mu = \int v\, d\mu >0$, so we can take 
the approximating functions $u_n$ and $v_n$ to be such that 
$\int u_n\, d\mu > 0$ and $\int v_n\, d\mu > 0$.  
The theorem produces needles $\nu_n$ of $\mu$ supported on $F$ and such that 
$$
\int u_n\, d\nu_n > 0, \qquad \int v_n\, d\nu_n > 0.
$$
Since $u \geq u_n$ and $v \geq v_n$, every such measure $\nu_n$ will be the 
required needle.

Let us now turn to the construction procedure.

Given $3$ affinely independent points $x,y,z$ in $E$, consider the linear continuous 
functionals $\lambda_x$ and $\lambda_y$ on $E$ introduced before Lemma 4.4
via the relations (4.1)-(4.2) and the Hahn-Banach theorem.
To each point $\theta \in \mathbb{S}^1 = \{(t,s): t^2 + s^2 = 1\}$ 
we can associate a linear functional
$\Lambda_\theta = t \lambda_x + s \lambda_y$ and define the function
$$
\Psi:\mathbb{S}^1 \to \mathbb{R}, \qquad
\theta = (t,s) \mapsto \int_{\{\ell_\theta(\xi - z) \geq 0\}} u(\xi)\, d\mu(\xi).
$$
Since $\mu\{\xi:\Lambda_\theta(\xi-z) = 0\} = 0$ (cf. (4.4)), this function is continuous 
on $\mathbb{S}^1$. In addition, we have the identity 
$$
\Psi(-\theta)+\Psi(\theta) = \int_E u\, d\mu. 
$$
Hence, by the intermediate value theorem, there exists $\theta \in \mathbb{S}^1$ 
such that
$$
\int_{\{\Lambda_\theta(\xi-z) \geq 0\}} u(\xi)\, d\mu(\xi) \, = \,
\int_{\{\Lambda_\theta(\xi-z) \leq 0\}} u(\xi)\, d\mu(\xi) \, = \,
\frac{1}{2}\,\int_E u\, d\mu.
$$
Also,
$$
\int_E v\,d\mu \, = \,
\int_{\{\Lambda_\theta(\xi-z) \geq 0\}} v(\xi)\, d\mu(\xi) + 
\int_{\{\Lambda_\theta(\xi-z) \leq 0\}} v(\xi)\, d\mu(\xi) \, > \, 0,
$$
so that at least one the last two integrals is positive.
Let $H^+$ denote one of the hyperspaces $\{\Lambda_\theta(\xi-z) \geq 0\}$ or
$\{\Lambda_\theta(\xi-z) \leq 0\}$ such that $\int_{H^+} v\, d \mu > 0$.
Necessarily, $\mu(H^+)>0$, and we may consider the normalized restriction
$\mu^+$ of $\mu$ to $H^+$ and will have the property that
\be
\int_{H^+} u\, d \mu^+ > 0, \quad
\int_{H^+} v\, d \mu^+ > 0.
\en

This procedure can be performed step by step along a sequence 
$\{(x_n,y_n,z_n)\}_{n \geq 1}$ of affinely independent points, chosen 
to be dense in $K \times K \times K$.
Let $\nu_1 = \mu^+$ be constructed according to the above procedure for 
$(x_1,y_1,z_1)$ and with an associated point $\theta_1 = (t_1,s_1) \in \mathbb{S}^1$.
Similarly, on the $n$-th step, given $\nu_n$, let $\nu_{n+1} = \nu_n^+$ be 
constructed for the triple $(x_n,y_n,z_n)$ and with the associated
linear functional 
$$
\Lambda_{\theta_n} = \Lambda_{(t_n,s_n)} = t_n \lambda_{x_n} + s_n \lambda_{y_n}.
$$

Since the space of all Borel probability measures on $K$ is compact and metrizable
for the weak topology, the sequence $\nu_n$ has a sub-sequential weak limit $\nu$. 
In particular, from (4.5) we derive the desired property
$$
\int_E u\, d\nu \geq 0,\qquad \int_E v\, d\nu \geq 0.
$$

It remains to show that ${\rm dim}(H_\nu)\leq 1$. Suppose not, in this case there 
exists affinely independent $x,y,z$ in the relative interior of $H_\nu$ that also
contains the points $2z - x$ and $2z - y$. Without loss of generality, let $z=0$, 
so that $\pm x$ and $\pm y$ belong to the relative interior of $H_\nu$. 
By the density property, there exists a subsequence, say $(x_k,y_k,z_k)$ such that 
$(x_k,y_k,z_k) \to (x,y,z)$. 

By the construction, the measure $\nu_k^+$ is supported on the half-space $H_k^+$, 
which is either $\{\xi:\Lambda_{(t_k,s_k)}(\xi-z_k) \geq 0\}$ or
$\{\xi:\Lambda_{(t_k,s_k)}(\xi-z_k) \leq 0\}$.
For definiteness, let it be the first half-space. Since all $H_k^+$ contain $x$ 
and $-x$, we then have
\be
\Lambda_{(t_k,s_k)}(x-z_k) \geq 0, \qquad \Lambda_{(t_k,s_k)}(-x-z_k) \geq 0,
\en
and similarly for the point $y$.

Recall that by Lemma 4.4, we can obtain a uniform bound $M$ such that 
$$
\|\Lambda_{(t_k,s_k)}\|_* \leq \|\lambda_{x_k}\|_* + \|\lambda_{y_k}\|_* \leq
M \qquad {\rm for \ all} \ k.
$$
Hence, $\Lambda_{(t_k,s_k)}(z_k) \rightarrow 0$ and 
$\Lambda_{(t_k,s_k)}(x_k - x) \rightarrow 0$ as $k \rightarrow \infty$. But then 
by (4.6), necessarily $\Lambda_{(t_k,s_k)}(x_k) \rightarrow 0$, as well. 
By the same argument, $\Lambda_{(t_k,s_k)}(y_k) \rightarrow 0$.

On the other hand, according to the definition of $\Lambda_{(t_k,s_k)}$ 
via (4.1)-(4.2), for each $k$,
$$
\Lambda_{(t_k,s_k)}(x_k-z_k) = t_k, \qquad \Lambda_{(t_k,s_k)}(y_k-z_k) = s_k,
$$ 
thus implying that $\lim_{k \to \infty}\, t_k = \lim_{k \to \infty}\, s_k = 0$. 
But this is impossible since $t_k^2 + s_k^2 = 1$.
This proves that ${\rm dim}(H_\nu)\leq 1$.
\qed


\vskip2mm
\section{Dilation and its properties}
\setcounter{equation}{0}

\vskip2mm
\noindent
Before turning to Theorem 1.3, we first comment on the basic properties of
the operation $A \rightarrow A_\delta$, where $\delta \in [0,1]$ 
is viewed as parameter. 

Let $F$ be a closed convex subset of a locally convex space 
$E$ with respect to which this operation is defined:
$$
A_\delta = \Big\{x \in A: m_\Delta(A) \geq 1 - \delta, 
\ {\rm for \ any \ interval} \ \Delta \subset F \ {\rm such \ that} 
\ x \in \Delta\Big\}.
$$
As before, $m_\Delta$ denotes a uniform distribution on $\Delta$ (understood 
as the Dirac measure, when the endpoints coincide). In this definiton, by the 
intervals $\Delta$ we mean closed intervals $[a,b]$ connecting arbitrary points 
$a,b$ in $F$. Moreover, the requirement $x \in \Delta$ may equivalently 
be replaced by the condition that $x$ is one of the endpoints of $\Delta$. 

Note that $A_1 = A$. If $0 \leq \delta < 1$, as an equivalent definition one could put
$$
A_\delta = \Big\{x \in F: m_\Delta(A) \geq 1 - \delta, 
\ {\rm for \ any \ interval} \ \Delta \subset F \ {\rm such \ that} 
\ x \in \Delta\Big\}.
$$
Indeed, in this case, if $x \in F \setminus A$, then $m_{[x,x]}(A) = 0 < 1 - \delta$ 
meaning that $x \notin A_\delta$ according to the second definition. Thus, 
for $\delta \in [0,1)$, both definitions lead to the same set and we have 
the property $A_\delta \subset A$.

\vskip5mm
{\bf Lemma 5.1.} {\it $a)$ 
If $A \subset F$ is closed, then every set $A_\delta$ is closed as well.

$b)$ If $E$ is a Fr\'echet space and $A$ is Borel measurable in $F$, 
then every set $A_\delta$ is universally measurable.
}

\vskip5mm
Let us recall that a set in a Hausdorff topological space $E$ is called 
universally measurable, if it belongs to the Lebesgue completion of the Borel
$\sigma$-algebra with respect to an arbitrary Borel probability measure on $E$.
In that case one may freely speak about the measures of these sets.

\vskip5mm
{\bf Proof.} For a Borel set $A$ in $F$, consider the function
$$
\psi(x,y) = \int 1_A\,dm_{[x,y]} = \int_0^1 1_A((1-t)x + ty)\,dt, \qquad x,y \in F.
$$

First assume that $A$ is closed. Then, given a net $x_i \rightarrow x$, 
$y_i \rightarrow y$ in $F$ indexed by a semi-ordered set $I$, we have 
$$
\limsup_{i \in I} 1_A((1-t)x_i + t y_i) \leq 1_A((1-t)x + ty). 
$$
After integration this implies
$$
\limsup_{i \in I}\, \psi(x_i,y_i) \, \leq \, \psi(x,y).
$$
Indeed, the space $L^1[0,1]$ is separable, so the above relation is only to be
checked for increasing sequences $i = i_n$ in $I$. But in that case
one may apply the Lebesgue dominated convergence theorem.
This means that $\psi$ is upper semicontinuous on $F \times F$, and thus 
$A_\delta$ represents the intersection over all $y \in F$ of the closed sets
$\{x \in A: \psi(x,y) \geq 1-\delta\}$.

In part $b)$, assume that $E$ is a Fr\'echet space. If $A$ is Borel, then
the function $\psi$ is Borel measurable on $F \times F$, so the complement of
$A_\delta$ in $A$,
$$
A \setminus A_\delta = \{x \in A: \psi(x,y) < 1-\delta, \
{\rm for \ some} \ y \in A\},
$$
represents the $x$-projection of a Borel set in $E \times E$. 
But every Borel set in a Polish space is Souslin, and therefore both 
$A \setminus A_\delta$ and $A_\delta$ are universally measurable
(cf. [Bog], Corollary 6.6.7 and Theorem 7.4.1).
\qed

\vskip2mm
There is an opposite operation representing a certain dilation 
or enlargement of sets. Given a Borel measurable set 
$B \subset F$ and $\delta \in [0,1)$, define
\be
B^\delta\ = \bigcup_{m_\Delta(B) > \delta} \Delta\ = \,
\big\{x \in F: m_{[x,y]}(B) > \delta \ \ {\rm for \ some} \ y \in F\big\}.
\en
Here the union is running over all intervals $\Delta \subset F$ such that 
$m_\Delta(B) > \delta$. 

Note that $B^\delta$ contains $B$
(since all singletons in $B$ participate in the above union).

\vskip5mm
{\bf Lemma 5.2.} {\it For any $\delta \in [0,1)$ and any Borel set $B \subset F$, 
the complement $A = F \setminus B$ satisfies the dual relations
$$
F \setminus A_\delta = (F \setminus A)^\delta \qquad and
\qquad F \setminus B^\delta = (F \setminus B)_\delta.
$$
In particular, $B^\delta$ is open in $F$, once $B$ is open in $F$.
}

\vskip5mm
{\bf Proof.} 
Given $x \in F$, the property $x \notin A_\delta$ means that, 
for some interval $\Delta \subset F$ containing $x$, we have 
$m_\Delta(A) < 1-\delta$, that is, $m_\Delta(B) > \delta$ meaning that 
$\Delta \subset B^\delta$. Therefore, 
$x \notin A_\delta \Leftrightarrow x \in B^\delta$.
For the last assertion, it remains to recall Lemma 5.1.
\qed

\vskip5mm
{\bf Lemma 5.3.} {\it Let $F$ be a convex closed set in $E$, and let $T$ 
be a linear continuous map from $E$ to another locally convex space $E_1$. 
For any Borel set $C \subset T(F)$,
$$
\big(T^{-1}(C) \cap F\big)^\delta = T^{-1}\big(C^\delta\big) \cap F,
$$
where the operation $C \rightarrow C^\delta$ is understood with respect to the
image $T(F)$.
}

\vskip5mm
{\bf Proof.} For all $a,b \in F$, the map $T$ pushes forward
the unform measure $m_{[a,b]}$ to $m_{[Ta,Tb]}$. Therefore,
the pre-image $B = T^{-1}(C)$ has measure $m_{[a,b]}(B) = m_{[Ta,Tb]}(C)$, so
$$
\big(B \cap F\big)^\delta\ = \bigcup_{m_{[Ta,Tb]}(C) > \delta} [a,b]
\ = \bigcup_{m_{[x,y]}(C) > \delta} T^{-1}([x,y]) \cap F \, = \,
T^{-1}\big(C^\delta\big) \cap F.
$$
\qed

\vskip2mm
When $E_1 = \R^n$ and $C$ is a polytope, the dilated set $C^\delta$ is a polytope, 
as well. Hence, by Lemma 5.3, $(T^{-1}(C))^\delta$ represents
the intersection of finitely many half-spaces.


\vskip2mm
\section{The dual form and proof of Theorem 1.3}
\setcounter{equation}{0}

\vskip2mm
\noindent
Following [B-N], let us reformulate Theorem 1.3 in terms of dilated sets. Putting 
$B = F \setminus A$ and using Lemma 5.2, the inequalty
\be
\mu(A) \, \geq \, 
\big[\,\delta \mu^*(A_\delta)^\alpha + (1-\delta)\big]^{1/\alpha} \qquad (0<\delta<1)
\en
is solved as $\mu_*(B^\delta) \geq R_\delta^{(\alpha)}(\mu(B))$, where
\be
R_\delta^{(\alpha)}(p) = 1 -
\bigg[\frac{(1 - p)^\alpha - (1-\delta)}{\delta}\bigg]^{1/\alpha}.
\en
More precisely,
in case $\alpha < 0$, the above expression is well-defined and represents
a strictly concave, increasing function in $p \in [0,1]$. For 
$\alpha = 0$, it is understood in the limit sense as
$$
R_\delta^{(0)}(p) = 1 - (1 - p)^{1/\delta}, \qquad 0 \leq p \leq 1,
$$
which is also strictly concave and increasing. In case $0 < \alpha \leq 1$, 
$R^{(\alpha)}(p)$ is defined to be (6.2) on the interval
$0 \leq p \leq 1 - (1 - \delta)^{1/\alpha}$ (when the expression makes 
sense) and we should put $R^{(\alpha)}(p) = 1$ on the remaining subinterval 
of $[0,1]$. 

In all cases, $R_\delta^{(\alpha)}:[0,1] \rightarrow [0,1]$ represents 
a concave, continuous, non-decreasing function such that 
$R_\delta^{(\alpha)}(0) = 0$ and $R_\delta^{(\alpha)}(1) = 1$. Put 
$R_0^{(\alpha)}(p) = \lim_{\delta \downarrow 0} R_\delta^{(\alpha)}(p) = 1$
for $0 < p \leq 1$ and $R_0^{(\alpha)}(0) = 0$.

\vskip5mm
{\bf Theorem 6.1.} {\it Let $\mu$ be an $\alpha$-concave probability measure 
on a complete locally convex space $E$ supported on a convex closed set $F$
$(-\infty < \alpha \leq 1)$. For any Borel subset $B$ of $F$ 
and for all $\delta \in [0,1)$,
\be
\mu_*(B^\delta) \, \geq \, R_\delta^{(\alpha)}(\mu(B)).
\en
}

\vskip5mm
For example, on the real line $E = \R$ for the Lebesgue measure $\mu$ on the 
unit interval $F = [0,1]$, we have $\alpha = 1$, and (6.3) becomes
$$
\mu(B^\delta) \, \geq \, \min\bigg\{\frac{1}{\delta}\ \mu(B), 1\bigg\}.
$$
For the Cauchy measures $\mu$ on $\R^n$ and $\R^\infty$ (cf. Examples 2.5), 
we have $\alpha = -1$, and then (6.2)-(6.3) with $F = E$ yield
$$
\mu(B^\delta) \, \geq \, \frac{\mu(B)}{1 - (1-\delta)(1-\mu(B))}.
$$

Note that when $E$ is a Fr\'echet space and $B$ is Borel, $B^\delta$ is universally 
measurable, so there is no need to use the inner masure.

Let us comment on the extreme values of $\delta$ in (6.1) and (6.3). Since the sets
$B^\delta$ increase for decreasing $\delta$, (6.3) will hold for $\delta=0$ 
by continuity, as long as this inequality holds for all $0<\delta<1$. In this case, 
(6.3) with $\delta = 0$ tells as that $\mu(B) > 0 \Rightarrow \mu(B^0) = 1$. 
Equivalently, after the substitution $A = F \setminus B$ and using Lemma 5.2, we get 
$\mu(A_0) = 0$, that is,
$$
\mu\Big\{x \in F: m_\Delta(A) = 1, 
\ {\rm for \ any \ interval} \ \Delta \subset F \ {\rm such \ that} 
\ x \in \Delta\Big\} \, = \, 0,
$$
as long as $\mu(A) < 1$. This case is however excluded from the formulation of 
Theorem 1.3 by the assumption $\mu^*(A_\delta)>0$. Note also that in case 
$\delta = 1$, (6.1) holds automatically, since then $A_1 = A$.

Thus, both Theorem 1.3 and Theorem 6.1 do not loose generality by assuming that 
$0 < \delta < 1$ (and we do this below in this section).

\vskip2mm
{\bf Equivalence of Theorem 1.3 and Theorem 6.1.} It is straightforward for 
$\alpha \leq 0$. This case also includes the values $\mu^*(A_\delta)=0$ in (6.1), 
since then $\mu_*(B^\delta) = 1$ for $B = F \setminus A$ and thus both (6.1) and (6.3) 
are immediate. 

Consider the case $0 < \alpha \leq 1$. For the implication $(6.1) \Rightarrow (6.3)$,
let $p = \mu(B)$, $0 < p < 1$. If $\delta \geq \delta_p = 1 - (1-p)^\alpha$, 
the formula (6.2) should be applied and then (6.3) becomes
\be
\mu_*(B^\delta) \, \geq \, 1 -
\bigg[\frac{(1 - \mu(B))^\alpha - (1-\delta)}{\delta}\bigg]^{1/\alpha}.
\en
Here the right-hand side tends to 1 as $\delta \downarrow \delta_p$, so necessarily
$\mu_*(B^{\delta_p}) = 1$ and hence $\mu_*(B^{\delta}) = 1$ for all 
$0 \leq \delta < \delta_p$. Thus, without loss of generality,
(6.3) may be stated as (6.4) for the range  $\delta \geq \delta_p$. 
If $\mu_*(B^\delta) = 1$ there is nothing to prove. If $\mu_*(B^\delta) < 1$,
then $\mu^*(A_\delta) > 0$ for the set $A = F \setminus B$. In that case, (6.1)
is exactly the same as (6.4).

For the implication $(6.3) \Rightarrow (6.1)$, assume that $\mu^*(A_\delta) > 0$.
Then $\mu_*(B^\delta) < 1$ for $B = F \setminus A$ which implies that
$\mu(B) < p_0 = 1 - (1-\delta)^{1/\alpha}$ according to the definition of 
$R_\delta^{(\alpha)}(\mu(B))$. Moreover, again the formula (6.2) 
should be applied to rewrite the hypothesis (6.3) in the form (6.4),
which can in turn be rewritten as (6.1).
\qed

\vskip5mm
{\bf Proof of Theorem 6.1.} Using Theorem 1.2, let us show how to reduce
the desired statement (6.3) to dimension one. Since the sets $B^\delta$ may only become
larger, when $F$ is getting larger, one may assume that $F = H_\mu$, i.e.,
the support of $\mu$. Fix $0 < \delta < 1$.

\vskip2mm
{\it Step} 1: First suppose that $B$ is an open set in $F$
such that the boundary $\partial B^\delta$ of $B^\delta$ in $F$ has $\mu$-measure zero.
Fix an arbitrary $p \in (0,1)$. Using the continuity of the functions 
$R_\delta^{(\alpha)}$, it is sufficient to show that
$\mu(B) > p \Rightarrow \mu(D) \geq  R_\delta^{(\alpha)}(p)$, where $D$ is the
closure of $B^\delta$. If this were not true, we would have
$$
\int (1_B - p)\,d\mu > 0, \qquad 
\int (R_\delta^{(\alpha)}(p) - 1_D)\,d\mu > 0,
$$
which is exactly the condition (1.5) for $u = 1_B - p$ and 
$v = R_\delta^{(\alpha)}(p) - 1_D$ (where $1_A$ denotes the indicator function
of a set $A$). These functions are lower-semicontinuous, so we may apply Theorem 1.2: 
For some one dimensional $\alpha$-concave probability measure $\nu$ supported on 
an interval $\Delta \subset F$, we have (1.6), i.e., 
$$
\nu(B) > p, \qquad \nu(D) < R_\delta^{(\alpha)}(\nu(B)).
$$
But
$$
\nu(D) \geq \nu\big(B^\delta\big) \geq \nu\big((B \cap \Delta)^\delta\big),
$$
where $(B \cap \Delta)^\delta$ is the result of the one dimensional dilation 
operation applied to $B \cap \Delta$ with respect to $\Delta$. Hence, we obtain
$\nu((B \cap \Delta)^\delta) < \nu(B)$ which contradicts the relation (6.3) 
in dimension one.

\vskip2mm
{\it Step} 2: Here we describe one class of open sets to which the previous step
may be applied. Let $B$ be a set of the form $T^{-1}(C) \cap F$, where
$T:E \rightarrow \R^n$ is a continuous linear map and $C \subset \R^n$ is 
an open polytope ($n \geq 1$ is arbitrary). Then $(T^{-1}(C))^\delta$ represents
an intersection of finitely many open half-spaces (Lemma 5.3). 
If $\mu(B)>0$,
then, by the zero-one law, the boundaries of these half-spaces have 
$\mu$-measure zero and hence $\mu(\partial B^\delta) = 0$, as well.

More generally, let $B = T^{-1}(C) \cap F$, where $C$ is a finite union of 
open polytopes in $\R^n$. Then $C^\delta$ is also a finite union of open polytopes.
Using Lemma 5.3, we obtain that 
${\rm clos}(B^\delta) \subset T^{-1}({\rm clos}(C^\delta))$, 
so $\partial B^\delta \subset T^{-1}(\partial C^\delta)$. Again
$\partial C^\delta$ is contained in finitely many hyperplanes of $\R^n$ and thus
$\mu(\partial B^\delta) = 0$.

\vskip2mm
{\it Step} 3: $B$ is an arbitrary open set in $F$, assuming that $F$ is 
a convex, compact set. Denote by $\cal G$ the collection of all cylindrical sets 
in $F$ described on the last step. Such sets constitute a base in the original
topology on $F$, since the two coincides by the compactness assumption. 
Hence $B = \cup\, G$, where the union is over all $G \in \cal G$ such that 
$G \subset B$. Since $\cal G$ is closed under finite
unions, one may apply the Radon property which gives
$$
\mu(B) = \sup\{\mu(G): G \in {\cal G}, \ G \subset B\}.
$$
For any $G$ as above, we have $\mu(G^\delta) \, \geq \, R_\delta^{(\alpha)}(\mu(G))$,
by the previous steps. Hence, we obtain (6.3) for $B$, as well.

\vskip2mm
{\it Step} 4: $B$ is an arbitrary Borel set in $F$. By the strengthened 
Radon property (1.7), it is sufficient to consider the case of a non-empty compact 
set $B$, and we may additionally assume that $F$ is compact.

Any open set in $F$ containing $x \in B$ contains this point together with 
$B(x) \cap F$, where $B(x) = T_x^{-1}(C(x))$. Here
$T_x:E \rightarrow \R^n$ is a continuous linear map and $C(x)$
is a Euclidean ball in $\R^n$ (with some $n$ depending on $x$). 
Using compactness of $B$, one can compose its finite covering by the sets 
of the form 
$$
G = \big(B(x_1) \cup \dots \cup B(x_N)\big) \cap F, \qquad x_j \in B,
$$ 
with full intersection being $B$.
Let $\{U_i\}_{i \in I}$ be a decreasing net indexed by a semi-ordered directed set 
$I$ such that each $U_i$ represents the intersection of finitely many sets $G$ 
as above. The latter guarantees that $\mu(U_i) \downarrow \mu(B)$ along the net.

Now, let $\delta < \delta' < 1$. Given $x \in F$, the property $x \notin B^\delta$ 
means that $m_{[x,y]}(B) = \inf_{i \in I}\, m_{[x,y]}(U_i) \leq \delta$
for any $y \in F$. In that case, there is $i \in I$ 
such that $m_{[x,y]}(U_i) < \delta'$, and hence the increasing sets
$$
V_i(x) = \big\{y \in F: m_{[x,y]}(U_i) < \delta'\big\}, \quad i \in I,
$$
cover $F$. By the construction, for each $i$, the function 
$\varphi(y) = m_{[x,y]}(U_i)$ is of the type
\bee
\varphi(y) \, = \, 
{\rm mes}\Big\{t \in (0,1): \forall\,k \leq l \ \ \exists j \leq N_k \ \ \ 
(1-t)T_{x_{kj}}(x) + tT_{x_{kj}}(y) \in C(x_{kj})\Big\}
\ene
for some continuous linear maps $T_{x_{kj}}:E \rightarrow \R^{n(x_{kj})}$ 
and some Euclidean balls $C(x_{kj})$ in $\R^{n(x_{kj})}$. 
As the boundaries of Euclidean balls do not contain non-degenerate intervals,
any such function $\varphi$ must be continuous on $F$. Therefore, all the sets
$V_i(x)$ are open in $F$, so that by compactness, $V_i(x) = F$ for some 
$i = i(x)$. 
Thus, given $x \in F \setminus B^\delta$, we have $m_{[x,y]}(U_{i(x)}) < \delta'$ 
for any $y \in F$, and hence $F \setminus B^\delta$ is contained in
$$
\bigcup_i  \big\{x \in F: m_{[x,y]}(U_i) < \delta' \ \ {\rm for \ all} \ y \in F\big\}.
$$
It follows that $B^\delta$ contains the intersection of the open sets
$$
U_i^{\delta'} \ = \big\{x \in F: m_{[x,y]}(U_i) > \delta' \ \ 
{\rm for \ some} \ y \in F\big\}
$$
and thus, by the Radon property,
$$
\mu_*(B^\delta)\, \geq\, \mu\big(\cap_i U_i^{\delta'}\big) =
\, \lim_{i} \mu\big(U_i^{\delta'}\big).
$$

On the other hand, by Step 3, 
$\mu(U_i^{\delta'}) \geq R_{\delta'}^{(\alpha)}(\mu(U_i))$, and taking the limit
along the net we get $\mu_*(B^\delta) \geq R_{\delta'}^{(\alpha)}(\mu(B))$.
It remains to let $\delta' \downarrow \delta$ and use the
contunuity of $R_{\delta}^{(\alpha)}$ with respect to $\delta$.
\qed


\vskip2mm
\section{Large and small deviations}
\setcounter{equation}{0}

\vskip2mm
\noindent
As is known, the dilation-type inequality (1.8) of Theorem 1.3 may equivalently be 
stated on functions (which is often more convenient in applications). 
Namely, with every Borel measurable function $u$ on $E$ with values in the
extended line $[-\infty,\infty]$, one associates its "modulus of regularity"
$$
\delta_u(\ep) \, = \, \sup\,
{\rm mes}\big\{t \in (0,1): |\,u((1-t)x + ty)\,| \leq \ep\, |u(x)|\big\}, 
\qquad 0 \leq \ep \leq 1,
$$
where the supremum is running over all points $x,y \in E$ such that $u(x)$ is finite.

The behavior of $\delta_u$ near zero is used to control the probabilities of large 
and small deviations of $u$ under hyperbolic measures by involving the parameter
$\alpha$, only (cf. [B4], [B-N], [F]). In particular, there is 
the following recursive functional inequality, which 
is stated below, in the setting of an abstract complete locally convex space $E$.

We assume that $\mu$ is an $\alpha$-concave probability measure on $E$ with 
$-\infty < \alpha \leq 1$ and that $u$ is a Borel measurable,
$\mu$-a.e. finite function on $E$.

\vskip5mm
{\bf Theorem 7.1.} {\it Given $0 < \lambda < {\rm ess\ sup}\, |u|$, for all 
$\ep \in (0,1)$,
\be
\mu\{|u| > \lambda \ep\} \, \geq \,
\Big[\delta\, \mu\{|u| \geq \lambda\}^\alpha + (1 - \delta)\Big]^{1/\alpha},
\en
where $\delta = \delta_u(\ep)$.
}

\vskip5mm
In case $\alpha = 0$, this relation turns into
\be
\mu\{|u| > \lambda \ep\} \, \geq \, \big(\mu\{|u| \geq \lambda\}\big)^\delta.
\en
Note that for $\alpha \leq 0$, the assumption 
$\lambda < {\rm ess\, sup}\, |u|$ may be removed. 

If $\mu$ is supported on a convex closed set $F$ in $E$, the inequalities (7.1)-(7.2)
continue to hold when $u$ is defined on $F$ (rather than on the 
whole space). In that case, in the definition of $\delta_u$ the supremum should 
be taken over all points $x,y \in F$.

\vskip5mm
{\bf Proof of Theorem 7.1.} Let us recall a simple argument based on Theorem 1.3. 
The latter is applied with $F = E$ to the set
$$
A = \{x \in E: \lambda \ep < |u(x)| < \infty\}.
$$
By the definiton,
\bee
A_\delta
 & = &
\{x \in E: m_{[x,y]}(A) \geq 1-\delta \ \ \forall y \in E\} \\
 & = &
\big\{x \in E: {\rm mes}\{t \in (0,1): \lambda \ep < |u((1-t)x + ty)| < \infty \} \geq 
1-\delta \ \ \forall y \in E\big\}.
\ene
Suppose that $\lambda \leq |u(x)| < \infty$. Then, for any $y \in E$, we have
$
|u((1-t)x + ty)| \leq \lambda \ep \, \Rightarrow \, |u((1-t)x + ty)| \leq \ep |u(x)|,
$
so that
$$
\hskip-35mm{\rm mes}\big\{t \in (0,1):|u((1-t)x + ty)| \leq \lambda \ep\big\}\ \leq
$$
$$
\hskip20mm
{\rm mes}\big\{t \in (0,1):|u((1-t)x + ty)| \leq \ep |u(x)|\big\}\, \leq\, \delta_u(\ep).
$$
Hence,
$$
{\rm mes}\big\{t \in (0,1): \lambda \ep < |u((1-t)x + ty)| < \infty\big\}\, \geq\,
1 - \delta_u(\ep)
$$
which implies that $x \in A_\delta$ with $\delta = \delta_u(\ep)$. 
This gives the inclusion 
$$
\{x \in E: \lambda \leq |u(x)| < \infty\} \subset A_\delta
$$ 
and also that $\mu_*(A_\delta) > 0$ (due to the assumption on $\lambda$).
It remains to apply (1.8).
\qed

\vskip2mm
In the next two corrolaries we follow [B-N], cf. also [F].
Denote by $m$, a median of $|u|$ under $\mu$, i.e., a real number such that
$$
\mu\{|u| > m\} \leq \frac{1}{2}, \quad \mu\{|u| < m\} \leq \frac{1}{2}.
$$

\vskip5mm
{\bf Corollary 7.2.} {\it Assuming that $m>0$, for all $r > 1$,
\be
\mu\{|u| \geq mr\} \, \leq
\bigg[1 + \frac{2^{-\alpha} - 1}{\delta_u(\frac{1}{r})}\bigg]^{1/\alpha}.
\en
}

\vskip5mm
When $\alpha = 0$, the right-hand side is understood as the limit
at zero, that is,
\be
\mu\{|u| \geq mr\} \, \leq\, 2^{-1/\delta_u(\frac{1}{r})}.
\en

If $\alpha < 0$, the inequality (7.3) may be simplified as
\be
\mu\{|u| \geq m r\} \, \leq\, C_\alpha\, \delta_u(1/r)^{-1/\alpha}
\en
with constant $C_\alpha = (2^{-\alpha} - 1)^{1/\alpha}$. 
Note $C_\alpha \rightarrow \frac{1}{2}$, as $\alpha \rightarrow -\infty$.
As is easy to see, we also have a uniform bound, such as, 
for example, $C_\alpha \leq 1$ in the region $\alpha \leq -1$.

\vskip5mm
{\bf Proof.} To derive (7.3) in case $\alpha \neq 0$, apply (7.1) with 
$\lambda = mr$ and $\ep = 1/r$. 
Then $\mu\{|u| > \lambda \ep\} \leq \frac{1}{2}$, and letting 
$p = \mu\{|u| \geq \lambda\}$, we get
$\frac{1}{2} \geq (\delta p^\alpha + (1 - \delta))^{1/\alpha}$.
It remains to solve this inequality in terms of $p$. Note that when $\alpha > 0$, 
necessarily $\frac{1}{2} \geq (1 - \delta)^{1/\alpha}$ 
or $\frac{2^{-\alpha} - 1}{\delta} \geq -1$,
so the right-hand side of (7.3) makes sense. By a similar argument,
(7.4) follows from (7.2) in the log-concave case.
\qed

\vskip2mm
{\bf Remark.}
An inequality of the form (7.5) can also be obtained by using a transportation 
argument, cf. [B5]. With this argument, a slightly weaker variant of (7.4) is
derived in [B4].

\vskip2mm
Now, let us turn to the problem of small deviations.

\vskip5mm
{\bf Corollary 7.3.} {\it If $m>0$, for all $0 < \ep < 1$,
\be
\mu\{|u| \leq m \ep\} \, \leq \, C_\alpha\, \delta_u(\ep)
\en
with constant $C_\alpha = \frac{2^{-\alpha} - 1}{-\alpha}$.
}

\vskip5mm
{\bf Proof.} One may assume that $\alpha \neq 0$ and $m=1$. From (7.1) with 
$\lambda = 1$, we obtain that $\mu\{|u| \leq \ep\} \leq  \phi(x)$, where
$\phi(x) = 1 - (1+x)^{1/\alpha}$ and $x = (2^{-\alpha} - 1)\,\delta_u(\ep)$.
Since this function is concave in $x>-1$, we have
$\phi(x) \leq \phi(0) + \phi'(0)x = \frac{2^{-\alpha} - 1}{-\alpha}\ \delta_u(\ep)$.
When $\alpha = 0$, (7.6) holds with 
$C_0 = \lim_{\alpha \rightarrow 0} C_\alpha = \log 2$.
\qed

\vskip5mm
Finally, let us illustrate Corollaries 7.2-7.3 on the example of the semi-norms.

\vskip5mm
{\bf Lemma 7.4.} {\it If $u$ is a Borel measurable semi-norm on $E$ $($not identically 
zero$)$, then 
$$
\delta_u(\ep) = \frac{2\ep}{1 + \ep}, \qquad 0 < \ep \leq 1.
$$
}

{\bf Proof.} One may assume that both $u(x)$ and $u(y)$ are finite in the definition
of $\delta_u$. Moreover, it is a matter of normalization alone, to assume that
$c = u(y) \leq u(x) = 1$. Then, by the triangle inequality,
$$
u((1-t)x + ty) \, \geq \, |(1-t) u(x) - tu(y)| \, = \, |(1+c)t - 1|, 
$$
so
\bee
{\rm mes}\big\{t \in (0,1): u((1-t)x + ty) \leq \ep\,u(x)\big\}
 & \leq &
{\rm mes}\big\{t \in (0,1):|(1+c)t - 1| \leq \ep \big\} \\
 & = & \min\{t_1,1\} - t_0,
\ene
where $t_1 = \frac{1+\ep}{1+c}$, $t_0 = \frac{1-\ep}{1+c}$. In case $c \geq \ep$,
we have $t_1 - t_0 = \frac{2\ep}{1+c} \leq \frac{2\ep}{1+\ep}$. In case $c \leq \ep$,
similarly $1 - t_0 = \frac{c+\ep}{1+c} \leq \frac{2\ep}{1+\ep}$. Thus, 
$\delta_u(\ep) \leq \frac{2\ep}{1+\ep}$ in both cases. Here, the equality is attained
by taking $y = -x$ with $0 < u(x) < \infty$.
\qed

\vskip5mm
Any Borel measurable semi-norm $u$ on $E$ is generated by a centrally symmetric, 
Borel measurable, convex set $B$ in $E$, so that
$$
B = \{x \in E: u(x) \leq 1\}.
$$
Moreover, $u$ is $\mu$-a.e. finite, if and only if $\mu(B)>0$ in which case the linear
hull of $B$ has $\mu$-measure 1 (by the zero-one law). We are then in position
to apply Corollary 7.2. More conveniently, starting from (7.1) with
$\lambda = r$ and $\ep = \frac{1}{r}$ ($r > 1$), Lemma 7.4 gives
\bee
1 - \mu(B)
 & = &
\mu\{u(x) > 1\} \\
 & \geq &
\big[\delta\, \mu\{u(x) \geq r\}^\alpha + (1 - \delta)\big]^{1/\alpha} \\
 & \geq &
\big[\delta\, (1 - \mu(rB))^\alpha + (1 - \delta)\big]^{1/\alpha}, \qquad
\delta = \frac{2}{r + 1}.
\ene
At this step, the assumption $\mu(B)>0$ may be removed. Recalling also Corollary 7.3,
we arrive at:

\vskip5mm
{\bf Corollary 7.5.} {\it Given a symmetric, Borel measurable, convex 
set $B$ in $E$, for all $r>1$,
\be
1 - \mu(B)  \, \geq \,
\bigg[\,\frac{2}{r+1}\ \big(1 - \mu(rB)\big)^\alpha + \frac{r-1}{r+1}\,\bigg]^{1/\alpha}.
\en
}

\vskip2mm
In the limit case $\alpha = 0$, the above is the same as
$$
1 - \mu(rB) \, \leq \, \big(1 - \mu(B)\big)^{(r+1)/2}.
$$
This inequality is due to Lov\'asz and Simonovits [L-S] in case of Euclidean balls 
$B$ in $\R^n$. Gu\'edon [G] extended it to general symmetric convex sets in $\R^n$ 
and also found a precise relation in the case $\alpha > 0$. 
Namely, (7.7) is solved in terms of $1 - \mu(rB)$ as
$$
1 - \mu(rB)\, \leq\, {\max}^{1/\alpha}\bigg\{\frac{r+1}{2}\ 
\big(1 - \mu(B)\big)^\alpha - \frac{r-1}{2},\, 0 \bigg\}.
$$
As for the range $\alpha < 0$, (7.7) may be then rewritten as
$$
1 - \mu(rB)\, \leq\, 
\bigg[\,\frac{r+1}{2}\ \big(1 - \mu(B)\big)^\alpha - \frac{r-1}{2}\,\bigg]^{1/\alpha}.
$$

These large deviations bounds provide a sharp form of Borell's Lemma 3.1 in
[Bor1]. 

Let us also mention an immediate consequence from Corollary 7.3 and 
Lemma 7.4 concerning measures of small balls.

\vskip5mm
{\bf Corollary 7.6.} {\it Given a symmetric, Borel measurable, convex 
set $B$ in $E$ such that $\mu(B) \leq \frac{1}{2}$, we have
$$
\mu(\ep B)  \, \leq \, C_\alpha\,\ep \qquad (0 \leq \ep \leq 1)
$$
with constant $C_\alpha = \frac{2(2^{-\alpha} - 1)}{-\alpha}$.
}

\vskip10mm
{\bf Acknowledgement.} We would like to thank V. Bogachev for careful reading
of the manuscript and valuable comments.


\vskip5mm
\begin{thebibliography}{Ma2}


\itemsep=-0.5pt
\small
\bibitem[B1]{B1}
Bobkov, S. G. Remarks on the growth of $L^p$-norms of polynomials. Geometric aspects 
          of functional analysis, 27–-35, Lecture Notes in Math., 1745, Springer, 
          Berlin, 2000. 

\bibitem[B2]{B2}
Bobkov, S. G. Some generalizations of Prokhorov's results on Khinchin-type inequalities 
          for polynomials. (Russian) Teor. Veroyatnost. i Primenen. 45 (2000), no. 4, 
          745--748. Translation in: Theory Probab. Appl. 45 (2002), no. 4, 644–-647. 

\bibitem[B3]{B3}
Bobkov, S. G. Localization proof of the isoperimetric Bakry-Ledoux inequality and 
          some applications. Teor. Veroyatnost. i Primenen. 47 (2002), no. 2, 340--346. 
          Translation in: Theory Probab. Appl. 47 (2003), no. 2, 308–-314. 

\bibitem[B4]{B4}
Bobkov, S. G. Large deviations via transference plans. Advances in mathematics research, 
          Vol. 2, 151–-175, Adv. Math. Res., 2, Nova Sci. Publ., Hauppauge, NY, 2003. 

\bibitem[B5]{B5}
Bobkov, S. G. Large deviations and isoperimetry over convex probability 
          measures. Electron. J. Probab. 12 (2007), 1072--1100.

\bibitem[B6]{B6}
Bobkov, S. G. On isoperimetric constants for log-concave probability distributions. 
          Geometric aspects of functional analysis, 81–-88, Lecture Notes in Math., 
          1910, Springer, Berlin, 2007. 

\bibitem[B-M]{B-M}
Bobkov, S. G., Madiman, M. Concentration of the information in data with log-concave 
          distributions. Ann. Probab. 39 (2011), no. 4, 1528–-1543. 

\bibitem[B-N]{B-N}
Bobkov, S. G., Nazarov, F. L. Sharp dilation-type inequalities with fixed 
          parameter of convexity. J. Math. Sci. (N.Y.) 152 (2008), no. 6, 826–-839.
          Translation from: Zap. Nauchn. Sem. POMI 351 (2007), Veroyatnost i 
          Statistika, 12, 54--78. 
           
\bibitem[Bog]{Bog}
Bogachev, V. I. Measure Theory. Vol. I, II. Springer-Verlag, Berlin, 2007. 
          Vol. I: xviii+500 pp., Vol. II: xiv+575 pp. 

\bibitem[B-S-S]{B-S-S} 
Bogachev, V. I., Smolyanov, O. G., Sobolev, V. I. Topological vector spaces and their
         applications (Russian). Moscow, Izhevsk, 2012, 584 pp.

\bibitem[Bor1]{Bor1}
Borell, C. Convex measures on locally convex spaces.
          Ark. Math. 12 (1974), 239--252.

\bibitem[Bor2]{Bor2} 
Borell, C. Convex set functions in $d$-space. Period. Math. Hungar.
          6 (1975), no. 2, 111--136.

\bibitem[B-L]{B-L} 
Brascamp, H. J., Lieb, E. H. On extensions of the Brunn-Minkowski and 
          Pr\'ekopa-Leindler theorems, including inequalities for log concave 
          functions, and with an application to the diffusion equation. 
          J. Funct. Anal. 22 (1976), no. 4, 366–-389. 

\bibitem[B-Z]{B-Z} 
Burago Yu. D., Zalgaller, V. A. Geometric inequalities. Springer-Verlag, Berlin, 
          1988. Translated from the Russian by A. B. Sosinskii, Springer Series 
          in Soviet Mathematics, xiv+331 pp.

\bibitem[D-K-H]{D-H-K} 
Davidovich, Ju. S., Korenbljum, B. I., Hacet, B. I. A certain property of 
          logarithmically concave functions. (Russian) Dokl. Akad. Nauk SSSR 
          185 (1969), 1215–-1218. 

\bibitem[F]{F}
Fradelizi, M. Concentration inequalities for $s$-concave measures of dilations of 
          Borel sets and applications. Electron. J. Probab. 14 (2009), no. 71, 
          2068–-2090.

\bibitem[F-G1]{F-G1} 
Fradelizi, M., Gu\'edon, O. The extreme points of subsets of s -concave probabilities 
          and a geometric localization theorem. Discrete Comput. Geom. 31 (2004), 
          no. 2, 327–-335.

\bibitem[F-G2]{F-G2} 
Fradelizi, M., Gu\'edon, O. A generalized localization theorem and geometric 
          inequalities for convex bodies. Adv. Math. 204 (2006), no. 2, 509–-529. 

\bibitem[G-M]{G-M} 
Gromov, M., Milman, V. D. Generalization of the spherical isoperimetric inequality 
          to uniformly convex Banach spaces. Composition Math. 62 (1987), 263--282.

\bibitem[G]{G} 
Gu\'edon, O. Kahane-Khinchine type inequalities for negative exponent. 
          Mathematika 46 (1999), no. 1, 165–-173. 

\bibitem[H-O]{H-O}
Hadwiger, H., Ohmann, D. Brunn-Minkowskischer Satz und Isoperimetrie. Math. Z., 
          66 (1956), 1--8.

\bibitem[I]{I} 
Ibragimov, I. A. On the composition of unimodal distributions. (Russian) Teor. 
          Veroyatnost. i Primenen. 1 (1956), 283–-288.

\bibitem[K-L-S]{K-L-S}
Kannan, R., Lov\'asz, L. Simonovits, M. Isoperimetric problems for convex bodies and 
          a localization lemma. Discrete Comput. Geom. 13 (1995), no. 3--4, 541–-559.

\bibitem[K-A]{K-A} 
Kantorovich, L. V., Akilov, G. P. Functional Analysis. Translated from the Russian 
          by Howard L. Silcock. Second edition. Pergamon Press, Oxford-Elmsford, 
          N.Y., 1982. xiv+589 pp. 

\bibitem[L-S]{L-S}
Lov\'asz, L. Simonovits, M. Random walks in a convex body and an improved volume 
          algorithm. Random Structures Algor. 4 (1993), no. 4, 359–-412.

\bibitem[M]{M}
Meyer, P.-A. Probability and potentials. Blaisdell Publishing Co. Ginn and Co., Waltham, 
          Mass.-Toronto, Ont.-London, 1966 xiii+266 pp. 

\bibitem[N-S-V]{N-S-V}
Nazarov, F., Sodin, M., Vol'berg, A. The geometric Kannan-Lov\'asz-Simonovits lemma, 
          dimension-free estimates for the distribution of the values of polynomials, 
          and the distribution of the zeros of random analytic functions. (Russian) 
          Algebra i Analiz 14 (2002), no. 2, 214--234. Translation in: 
          St. Petersburg Math. J. 14 (2003), no. 2, 351–-366. 

\bibitem[P-W]{P-W}
Payne, L. E., Weinberger, H. F. An optimal Poincar\'e inequality for convex domains. 
          Arch. Rational Mech. Anal. 5 (1960), 286–-292. 

\bibitem[Ph]{Ph}
Phelps, R. R. Lectures on Choquet's theorem. D. Van Nostrand Co., Inc., Princeton, 
          N.J.-Toronto, Ont.-London, 1966 v+130 pp. 

\bibitem[Pr]{Pr}
Pr\'ekopa, A. Logarithmic concave measures with application to stochastic programming. 
          Acta Sci. Math. (Szeged) 32 (1971), 301–-316. 

\bibitem[R]{R}
Rudin, W. Functional Analysis. McGraw-Hill Series in Higher Mathematics. 
          McGraw-Hill Book Co., New York-D\"usseldorf-Johannesburg, 1973, xiii+397 pp. 

\end{thebibliography}
\end{document}